\documentclass[12pt]{amsart}
\usepackage{amsmath}\usepackage{amssymb}
\usepackage{amscd}
\usepackage[all]{xy}
\def\lra{\longrightarrow}
\def\map#1{\,{\buildrel #1 \over \lra}\,}
\def\smap#1{{\buildrel #1 \over \to}}

\newcommand{\pc}{{\mathrm{pc}}}
\newcommand{\iso}{{\mathrm{iso}}}
\newcommand{\into}{\rightarrowtail}
\newcommand{\onto}{\twoheadrightarrow}

\def\cA{\mathcal A}
\def\cB{\mathcal B}
\def\cC{\mathcal C}
\def\cE{\mathcal E}

\def\cM{\mathcal M}
\def\cS{\mathcal S}
\def\cT{\mathcal T}

\newcommand{\R}{\mathbb{R}}
\newcommand{\Q}{\mathbb{Q}}
\newcommand{\Z}{\mathbb{Z}}
\newcommand{\N}{\mathbb{N}}
\newcommand{\F}{\mathbb{F}}
\newcommand{\m}{\mathfrak{m}}
\newcommand{\Hom}{\operatorname{Hom}}
\newcommand{\Aut}{\operatorname{Aut}}
\newcommand{\bbS}{\mathbb{S}}

\newcommand{\Sets}{\mathbf{-Sets}}
\def\fin{{\text{f}}}

\input xy
\xyoption{all}

\numberwithin{equation}{section}

\theoremstyle{plain}
\newtheorem{thm}[equation]{Theorem}

\newtheorem{lem}[equation]{Lemma}
\newtheorem{cor}[equation]{Corollary}

\newtheorem{substuff}{Remark}[equation]

\theoremstyle{definition}
\newtheorem{defn}[equation]{Definition}
\newtheorem{ex}[equation]{Example}
\newtheorem{subex}[substuff]{Example}
\newtheorem{porism}[equation]{Porism}

\theoremstyle{remark}

\newtheorem{subrem}[substuff]{Remark}

\begin {document}
\title{The $K'$--theory of monoid sets}
\date{\today}
\author{Christian Haesemeyer}
\address{School of Mathematics and Statistics, University of Melbourne,
VIC 3010, Australia}
\thanks{Haesemeyer was supported by ARC DP-170102328}

\author{Charles A.\,Weibel}
\address{Math.\ Dept., Rutgers University, New Brunswick, NJ 08901, USA}
\email{weibel@math.rutgers.edu}\urladdr{http://math.rutgers.edu/~weibel}
\thanks{Weibel was supported by NSF grants}
\maketitle
{\vspace{2ex}%\begin{quotation}\begin{center}\begin{em}}
\begin{center}
\hspace{1ex}{C.\,A.\,Weibel dedicates this paper
to his father, C.\,G.\,Weibel.}
\end{center}
{\vspace{2ex}%{\par\end{em}\end{center}\end{quotation}}
\medskip

There are two flavors of $K$-theory for a pointed
monoid $A$ in the literature; see \cite{Deitmar}, \cite{CLS}, \cite{Mah} and \cite{Sz}
for example. One is $K(A)$, the $K$-theory of the
category of finitely generated projective $A$--sets;
the other is $G(A)$, the $K$-theory of the category of all
finitely generated $A$--sets. While $K(A)$ is quite simple,
$G(A)$ turns out to be complicated.

In this paper, we study $K'(A)$, the $K$-theory of the
category of finitely generated {\it partially cancellative}
$A$--sets (see Definition \ref{def:pc}). Such $A$--sets are
well-behaved, include free $A$--sets, and turn out to be useful
in toric geometry; see \cite{CHWW}.  An important example of a pc monoid is the free
pointed monoid $\N = \{\ast,1,t,t^2\!,...\}$.
Our key tools are the
additivity, devissage and localization theorems for ``CGW''-categories,
developed by Campbell and Zakharevich \cite{CZ}.

Our main result is an analogue of the ``Fundamental Theorem''
%\edit{we hope!}
in algebraic $K$-theory (Theorem \ref{Fund.Thm}): if $A$ is an abelian partially cancellative monoid then
\[
K'(A)\simeq K'(A\wedge\N) \quad \mathrm{and} \quad K'(A\wedge\Z_+)\simeq K'(A)\vee \Omega^{-1}K'(A).
\]
In particular, $K'(\N)\simeq \bbS$, the sphere spectrum.
The corresponding result fails dramatically for $G$-theory; see Example \ref{G-and-K'}.

%Several authors have studied the $K$-theory of the category of
%finitely generated projective sets, written as $K(A)$;
%see \cite{Deitmar}, \cite{CLS}.

We define partially cancellative monoids in Section \ref{sec:pc}.
In Section \ref{sec:K} we describe the $K$-theory of various
categories of $A$--sets, using the constructions and techniques of \cite{CZ}.
In Section \ref{sec:K'} we establish a number of structural theorems.
%localisation theorem \ref{A-A/s} for the $K'$-theory of monoids.
In Section \ref{sec:FT}  we prove the Fundamental
Theorem alluded to above, as well as a monoid version of the Farrell--Hsiang computation of the $K$-theory
of twisted polynomial rings \cite{FarrellHsiang}.  Finally, in Section \ref{sec:MS} we discuss the $K'$-theory of monoid schemes and in particular compute it for projective spaces.

\noindent{\it Notation:}
By a {\it monoid} $A$ we mean a pointed set with an associative
%commutative
product and distinguished elements $\ast,1$ such that
$a\cdot1=a = 1\cdot a$ and $a\cdot\ast=\ast = \ast\cdot a$ for all $a\in A$; these
are sometimes called ``pointed
%abelian
monoids.''
%We write $\N$ for the free pointed monoid $\{\ast,1,t,t^2,...\}$.
% and $\N/\!/n$ for the monoid
%$\N/(t^n\N)=\{1,t,...,t^{n-1}\}\amalg\{\ast\}$.
If $G$ is a group, we write $G_+$ for the pointed monoid $G\coprod\{\ast\}$. We note that the initial monoid %$\N/\!/1=
$\{\ast,1\}$ is called $\mathbb{F}_1$ in \cite{CLS}.

We shall always assume all categories are skeletally small. We shall write $\bbS$ for the sphere spectrum and $S^\infty X$ for the suspension spectrum of a space $X$.
%By abuse, we shall write $\Omega^\infty S^\infty$ for the sphere spectrum.

\smallskip

\section{Partially cancellative monoids}\label{sec:pc}
A (left) $A$--set is a pointed set $X$ with an action of $A$;
in particular, $1\cdot x=x$ and $\ast\cdot x=\ast$ for all $x\in X$.
For example, a free $A$--set is just a wedge of copies of $A$.

%Recall that a pointed abelian monoid $A$ is {\em cancellative}
%if whenever $ac = bc$ for some $c\ne\ast$, we have $a=b$.

%Equivalently, $B$ is cancellative if $B \setminus \{0\}$ is an
%unpointed monoid that maps injectively to its group completion
%$(B\setminus \{0\})^+$. In this situation, the
%{\em pointed group completion} of $B$ is
%$$
%B^+ := \left((B \setminus \{0\})^+\right)_*.
%$$

\begin{defn}\label{def:pc}
A monoid $A$ is {\em partially cancellative}, or {\em pc} for short,
if $A$ is noetherian and
$ac=bc\ne\ast$ (or $ca = cb \neq \ast$) implies $a=b$ for all $a,b,c$ in $A$.
The prototype finite pc monoid is $\{\ast,1,t,...,t^N=t^{N+1}\}.$
%\edit{pictures!}

If $A$ is a pc
%abelian
monoid, we say that a pointed (left) $A$--set $X$ is
{\em partially cancellative} if for every $x\in X$ and $a,b$ in $A$,
if $ax=bx\ne\ast$ then $a=b$.
Note that the subcategory of pc monoid sets is
closed under subobjects and quotients.
\end{defn}

\begin{subrem}\label{subrem:maximal}
	If $A$ is a pc monoid, then the subset $\mathfrak{m}$ of
	non-units in $A$ is a two-sided ideal, and is the unique maximal (two-sided) ideal of $A$.
	Indeed, if $xy = 1$, then $xyx = x$ and therefore $yx = 1$.
\end{subrem}

%A monoid $A$ is {\em partially cancellative and torsion free}, or
%{\em pctf} for short, if $A$ is isomorphic to $B/I$ where $B$ is a
%cancellative monoid whose group completion is torsion free and $I$
%is an ideal of $B$.

\begin{ex}
Let $\N$ denote the pointed monoid $\{\ast,1,t,t^2,...\}$.
A (pointed) $\N$--set is just a pointed set $X$ with a successor
function $x\mapsto tx$.
Thus we may identify a finite $\N$--set with
a (pointed) directed graph such that every vertex has outdegree~1.
Every finite rooted tree is a pc $\N$--set;
the successor of a vertex $x$ is the adjacent vertex closer to the root vertex $\ast$.
In fact, a finite $\N$--set is partially cancellative if and only
if it is a rooted tree, because for every $x\in X$, the sequence
$\{x,tx,t^2x,...\}$ terminates at the basepoint.

An $\N$--set is not pc if and only if it contains a loop, i.e.,
there is an element $x\ne\ast$ and an integer $d$ such that $t^dx=x$.
A typical non-pc $\N$--set is
$\{\ast,1,t,...,t^d,...,t^{N -1}\;\vert\; t^{N}=t^{d}\}$.
\end{ex}

\begin{ex}\label{G-sets}
Let $A$ be $G_+$, where $G$ is a group, with
a disjoint basepoint $\ast$ added. Then every $A$--set is a
wedge of copies of cosets $(G/H)_+$.
If $H$ is a proper subgroup of $G$, then $(G/H)_+$ is
not pc because $h\cdot H = 1\cdot H$ for $h\in H.$ %and $x = 1\,H.$
Thus an $A$--set is pc if and only if it is free.
%hence $K(A)=K'(A)$.
\end{ex}

\section{Quasi-exact categories}\label{sec:K}

\begin{defn}(\cite{Deitmar}, \cite[Ex.\,IV.6.14]{WK})
A {\it quasi-exact category} is a category $\cC$ with a distinguished zero object,
and a coproduct $\vee$, equipped with a family $\cS$ of sequences of the form
\begin{equation*}\label{admissible}
  0 \to Y \map{i} X  \map{j} Z \to 0,
\end{equation*}
called ``admissible,'' such that: (i) any sequence isomorphic to an
admissible sequence is admissible; (ii) for any admissible sequence,
$j$ is a cokernel for $i$ and $i$ is a kernel for $j$;
(iii) $\cS$ contains all split sequences (those with $X\cong Y\vee Z$);
and (iv) the class of admissible epimorphisms $j$ (resp.,
admissible monics $i$) is closed under composition and
pullback along admissible monics (resp., pullback along admissible
epimorphisms). We will often write $X/Y$ for the cokernel of $X\into Y$.

Quillen's construction in \cite{Q} yields a category $Q\cC$,
and $K(\cC)$ is the connective spectrum with initial space
$\Omega BQ\cC$; we write $K_n(\cC)$ for $\pi_nK(\cC)$.
The group $K_0(\cC)$ is generated by the objects of $\cC$, modulo the relations
that $[X]=[Y]+[Z]$ for every admissible sequence.
\end{defn}

\begin{subex}\label{BBQ}
The category $\mathbf{Sets}_\fin$ of finite pointed sets
is quasi-exact; every admissible sequence is split exact.
It is well known that the Barratt--Priddy--Quillen theorem implies that
$K(\mathbf{Sets}_\fin)\simeq  \bbS$
(see \cite{Deitmar}, \cite{CLS}, \cite[Ex.\,IV.6.15]{WK}.)

More generally, the category of finitely generated projective $A$--sets
is quasi-exact for any monoid $A$; every admissible exact sequence is split.
The $K$-theory of the category of finitely generated projective $A$--sets
is written as $K(A)$.
Following \cite[3.1]{Deitmar} and \cite[2.27]{CLS}, we see that if $A$ has no
idempotents or units then $K(A)\simeq \bbS$.
\end{subex}

\begin{subex}\label{G-and-K'}
If $A$ is a left noetherian monoid, the category
of finitely generated pointed left $A$--sets is quasi-exact; a sequence
\eqref{admissible} is admissible if $Y\into X$ is an injection, and
$Z$ is isomorphic to the quotient $A$--set $X/Y$.
(See \cite[Ex.\,IV.6.16]{WK}.)

Following \cite{CLS}, %Chu et al.\
we define $G(A)$ to be the $K$-theory of this category.
For example, the group $G_0(\N)=\pi_0G(\N)$ is the
free abelian group on $[\N]$ and the infinite set of loops of varying lengths.
%$\N$--sets $\N/t^n\N$.

If $A=G_+$ for a group $G$, then
$G_0(A)$ is the Burnside ring of $G$; for example,
$G_0(\Z_+)$ is free abelian on the classes of the (pointed)
cyclic groups $\Z_+$ and the loops $(\Z/n\Z)_+$.
By \cite[5.2]{CDD}, the spectrum $G(A)$ is the $G$--fixed points %spectrum
of the equivariant sphere spectrum, at least if $G$ is an abelian group.
%\edit{"abelian" because that's how their result is stated}
\end{subex}

\begin{subex}\label{K'-example}
If $A$ is a pc monoid, the category $A\Sets_\pc$ of
finitely generated pc $A$--sets is quasi-exact.  We write
$K'(A)$ for $K(A\Sets_\pc)$.
\end{subex}

%\begin{rem}
%$G_0(\N)$ is the free abelian group on $[\N]$ and the infinitely many
%\edit{move!}
%$\N$--sets $(\Z/n)_+$.
%\end{rem}

\begin{subex}\label{K(G)}
When $A=G_+$ for a group $G$, we saw in Example \ref{G-sets} that
every pc $A$--set is a free $A$--set. Therefore $K(G_+)=K'(G_+)$.

If $X=(G_+)^{\vee r}$, then $\Aut(X)$ is the wreath product $G\wr\Sigma_r$.
By the Barratt--Priddy--Quillen theorem,
$K(A)$ is the group completion of $\coprod B(G\wr\Sigma_r)$, which is
$\Omega^\infty S^\infty(BG_+)$. %by the Barratt--Priddy--Quillen theorem.
In particular when $A=\Z_+$, we have
\[
K'(\Z_+) \simeq K(\Z_+) \simeq \bbS\vee\Omega^{-1}\bbS, \quad K'_n(\Z_+) \cong \pi^s_n\oplus\pi^s_{n-1}
\]
This calculation of $K(G_+)$ is well known;  see \cite[5.9]{CLS}.%
\footnote{The formula in \cite[p.\,146]{Deitmar}
is incorrect, as $\Aut(G_+^{\vee n})$ is the wreath product $G\wr\Sigma_n$.}
\end{subex}

If $\cC$ is a quasi-exact category, we can form a double category
$(\cM,\cE)$ with the same objects as $\cC$; the horizontal and vertical
maps are the admissible monics and epis (composed backwards), respectively,
and the 2--cells are commutative diagrams of the form
\begin{equation}\label{ME-square}
\SelectTips{cm}{10 scaled 2000}\UseTips
\xymatrix@R=1.9em{
Y~ \ar@{>->}[r] & X \\
Y'~ \ar@{>->}[r]\ar@{->>}[u] & X'.\ar@{->>}[u]
}\end{equation}
We say that a square \eqref{ME-square}
is {\it distinguished} if the natural map of cokernels
$X'/Y'\to X/Y$ is an isomorphism.  Thus distinguished squares are both
pushout squares and pullback squares.
Note that $\cC$ is an ``ambient category'' for $(\cM,\cE)$ in the sense of
\cite[2.2]{CZ}.
We define
$k(X\onto X/Y)$ to be $(Y\into X)$ and $c(Y\into X)=(X\onto X/Y)$.

\begin{lem}
If $\cC$ is a quasi-exact category, $(\cM,\cE)$ is a CGW-category
in the sense of \cite[2.4]{CZ}.
\end{lem}

\begin{proof}
We need to verify the axioms in {\em loc.\,cit.}
The isomorphism $\iso(\cE)\cong\iso(\cM)$ is \cite[2.2]{CZ}.
Axiom (Z) holds because $\cC$ has a zero object, Axiom (I)
follows from (i) and (iii), and axiom (A) is (iii).
Axiom (K) is immediate from the definitions of $c$ and $k$, and
Axiom (M) is (ii).
\end{proof}

\begin{subrem}
In \cite{CZ}, Campbell and Zakharevich define the $K$-theory of a CGW-category using an appropriate version
of the $Q$-construction. If $\cC$ is a quasi-exact category with associated CGW-category $(\cM,\cE)$,
then $Q\cC$ and $Q(\cM,\cE)$ are isomorphic categories, and thus the two possible definitions of the $K$-theory
coincide.

With this in mind, we will abuse notation and consider a quasi-exact category as a CGW-category.
\end{subrem}

\begin{lem}\label{1st-isom}
If $X$ is a pointed $A$--set and $Y, Z$ are pointed $A$--subsets, then
$Y/(Y\cap Z) \into X/Z \onto X/(Y\cup Z)$ is an admissible sequence.
That is, $Y/(Y\cap Z) \smap{\cong} (Y\cup Z)/Z$.
\end{lem}

\begin{proof}
This is standard, and left to the reader.
\end{proof}

\begin{thm}\label{ACGW}
Let $A$ be a pointed monoid, and $\cC$ a quasi-exact subcategory of
$A\Sets$ closed under pushouts and pullbacks. Then the associated double
category $(\cM,\cE)$ is an ACGW-category in the sense of \cite[5.2--5.3]{CZ}.
\end{thm}

\begin{proof}
We let ``commutative square'' mean a 2-cell of the form \eqref{ME-square}.
For this, we need to check the axioms (P), (U), (S) and (PP).
Axiom (P) is the evident assertion that $\cM$ is closed under pullbacks,
and $\cE$ is closed under pushouts.

Given Lemma \ref{1st-isom},
Axiom (U) is equivalent to the following assertion:
given a pointed $A$--set $X$ and
pointed $A$--subsets $Y, Z$ (i.e., a diagram
$Y\into X\onto X/Z$), the pushout $(X/Z)\cup_X(X/Y)$ and pullback
$Y\times_X Z=Y\cap Z$ fit into a commutative diagram:
\[
\SelectTips{cm}{10 scaled 2000}\UseTips
\xymatrix@R=1.9em{
\square~ \ar@{>->}[r] & X/Z \ar@{->>}[r] & (X/Z)\cup_X(X/Y) \\
Y\mathstrut~ \ar@{->>}[u] \ar@{>->}[r] & X \ar@{->>}[u]\ar@{->>}[r]
& X/Y \ar@{->>}[u] \\
Y\cap Z~ \ar@{>->}[u] \ar@{>->}[r] & Z~ \ar@{>->}[u] \ar@{->>}[r]
& D\mathstrut \ar@{>->}[u]
}\]
(The object $\square$ is $Y/(Y\cap Z) \cong (Y\cup Z)/Z$.)
%by Lemma \ref{1st-isom}, there is an object $\square$ such that the
%top horizontal row and the left vertical row are admissible sequences,
%and the upper left square commutes.
%Thus (U) follows from the definition of commutative square.

Axiom (S) says that given a pullback
square in $\cM$ of the form
\begin{equation*}%\label{axiom-S}
\SelectTips{cm}{10 scaled 2000}\UseTips
\xymatrix@R=1.9em{
Y\cap Z~\mathstrut \ar@{>->}[r]\ar@{>->}[d] & Z\mathstrut\ar@{>->}[d] \\
Y~\  \ar@{>->}[r] & X,
}\end{equation*}
the pushout $P=Y\cup_{Y\cap Z}Z$ exists, and $X/P$ is the pushout of
$X/Y$ and $X/Z$ along $X/K$.  Dually, given a pushout square in $\cE$,
of $X\onto X/Y$ and $X\onto X/Z$, the pullback $L=X/(Y\cup Z)$
of $X/Y$ and $X/Z$ along $(X/Y)\cup_X(X/Z)$ exists, and the
kernel of $X\to L$ is $Y\cup Z$.  Both of these assertions are
elementary, in the spirit of Lemma \ref{1st-isom}.

Axiom (PP) says that if $V$ is an $A$--subset of both $X$ and $Y\!,$
then  $V$ is the intersection of $X$ and $Y$ in $X\cup_VY$.
(The pushout $X\cup_VY$ of $X$ and $Y$ along $V$ is the quotient
of the wedge $X\vee  Y\!,$ modulo the equivalence relation
identifying the two copies of $V$; it is the object
$X\!\star_V\!Y$ of \cite[5.3]{CZ}.)
%to be the pushout $X\cup_VY$ of $X$ and $Y$ along $V$,
%i.e., the wedge of $X$ and $Y$\!, modulo the equivalence relation
%identifying the two copies of $V$.  Then $V$ is the intersection
%of $X$ and $Y$ in $X\cup_VY$.
Dually, given epis $X\onto Z$ and $Y\onto Z$, the same argument
shows that the kernel of $(X\times_Z Y)\onto X$ is
isomorphic to the kernel of $Y\onto Z$.
\end{proof}

\begin{subrem}
It would suffice to prove Theorem \ref{ACGW} for pointed sets,
because the forgetful functor from $A\Sets$ to pointed sets
creates colimits and limits; its left adjoint is the free
functor, and its right adjoint sends a pointed set $X$ to its
co-induced $A$--set $\Hom(A,X)$. However, the proof would be no shorter.
\end{subrem}

\begin{cor}\label{cor-ACGW}
If $A$ is a left noetherian monoid, the associated double category $(\cM,\cE)$
associated to $A\Sets_{fg}$ (fin.\,gen.\,$A$--sets) is an ACGW--category.
The same is true for any subcategory closed under
pushouts and pullbacks, such as the categories
of pc $A$--sets (if $A$ is pc) and finite $A$--sets.
\end{cor}

\begin{proof}
It is straightforward that finitely generated (or pc, if $A$ is pc, or finite)
$A$--sets are closed under pushouts and pullbacks in the category $A\Sets$ when $A$ is
left noetherian. Thus Theorem \ref{ACGW} applies.
%If $V$ is an $A$--subset of both $X$ and $Y\!,$ we define
%$X\!\star_V\!Y$ to be the pushout $X\cup_VY$ of $X$ and $Y$ along $V$,
%i.e., the wedge of $X$ and $Y$\!, modulo the equivalence relation
%identifying the two copies of $V$.  Then $V$ is the intersection
%of $X$ and $Y$ in $X\cup_VY$.  Dually, given epis $X\onto Z$ and
%$Y\onto Z$, the same argument shows that the kernel of
%$X\times_VY\onto X$ is isomorphic to the kernel of
%\edit{dual st OK?}
%$Y\onto Z$.
\end{proof}

\section{$K'$-theory of a monoid}\label{sec:K'}

We recall the following definition from Example \ref{K'-example}.

\begin{defn}
If $A$ is a pc monoid, $K'(A)$ denotes the $K$-theory
of the category of finitely generated, partially cancellative $A$--sets.
\end{defn}

\begin{subrem}\label{K'-functor}
$K'$-theory is contravariantly functorial for normal monoid maps
(we recall from \cite{CLS} that a monoid map is {\it normal} if it is the
composition of a quotient by an ideal and an injection), and covariantly
functorial for flat monoid maps (a monoid map $A\to B$ is {\it flat}
if the base extension functor $X\mapsto B\wedge_AX$ is exact).
\end{subrem}

We will need the following theorem, taken from \cite[6.1]{CZ}.
\begin{thm}[Devissage]\label{thm:devissage}
%The Devissage Theorem \cite[6.1]{CZ} says that
Let $\cA$ be a full pre-ACGW category of a pre-ACGW category $\cB$,
closed under subobjects and quotient objects, and such that
the maps in $\cE_{\cA}$ are the maps in $\cE_{\cB}$ lying in $\cA$.
Suppose that every object of $\cB$ has a finite filtration
$$1\into F_1B\into\cdots\into F_nB=B$$
such that each $F_iB/F_{i-1}B$ is in $\cA$.
Then $K(\cA)\cong K(\cB)$.
\end{thm}

%\begin{lem}
%$K'(\N/\!/n) \simeq K'(\N/\!/1)\simeq \Omega^\infty S^\infty$ for all $n\ge1$.
%\end{lem}
%
%\begin{proof}
%Every pc $\N/\!/n$--set $X$ has a filtration by sub--sets
%$F_iX=\{x\in X | t^ix=\ast\}$, and the quotients $F_iX/F_{i+1}X$ are
%$\N/\!/1$--sets.  The lemma follows by Devissage.
%\end{proof}

\begin{lem}\label{lem:length}
If $A$ is a pointed pc monoid of finite length with units $G$, then
$$K'(A)\simeq S^\infty(BG_+).$$
In particular, if $A$ has no nontrivial units then
%If $A$ is a finite pointed monoid with no units then
$K'(A)\simeq\bbS$. In particular,
$K'(\N/t^n\N) \simeq\bbS$ for all $n\ge1$.
\end{lem}

\begin{proof}
The unique maximal ideal $\m$ of $A$ (see Remark \ref{subrem:maximal}) defines a finite filtration of any
finitely generated $A$--set $X$: if $\m^n=\ast$ then
\[
X \supset \m X \supset \m^2 X \supset\cdots \m^iX\supset\cdots\supset\m^nX=\ast.
\]
The lemma follows by Devissage \ref{thm:devissage}, because $A/\m=G_+$. %by hypothesis.
\end{proof}

The following Additivity Theorem is a special case of \cite[7.14]{CZ}.
Recall from Theorem \ref{ACGW} that the double category $(\cM,\cE)$
associated to a quasi-exact category is an ACGW-category
in many cases of interest.

\begin{thm}[Additivity]\label{additivity}
If $\cB$ is a quasi-exact category
such that $(\cM,\cE)$ is an ACGW-category,
and $s\into t \onto q$ is an admissible sequence of exact functors
$\cB \to \cA$, then %as functors $\cB\to\cA$:
\[ t_* = s_* + q_* : K(\cB) \to K(\cA). \]
\end{thm}
%should follow from \cite[7.14]{CZ}; $\cA$ and $\cB$ must be ACGW.
%Taking $\cA=\cB$, we have $K(S_1(\cB,\cA))\simeq K(\cB)\vee K(\cB)$.

\begin{proof}
By \cite[7.5, 7.11]{CZ}, there is a CGW-category $S_2\,\cB$
of admissible sequences $B_0\into B_1\onto B_2$.
The given admissible sequence defines exact functors
$s,q,t: S_2\,\cB\to\cA$, and if $\amalg:\cB \times \cB \to S_2\,\cB$ is the coproduct
functor then the following compositions agree:
\def\specialsqt{\underset{\scriptstyle s\vee q}
	{\overset{t}{\rightrightarrows}}}
\[
\cB \times \cB \map{\amalg} S_2\,\cB  \specialsqt  \cA.
\]
Hence they give the same map on $K$-theory.
By \cite[7.14]{CZ}, the source and target functors yield an equivalence
%\edit{end additivity}
\[ K(\cB)\vee K(\cB) \map{\simeq} K(S_2\cB). \]
Hence $t_* = (s\vee q)_* = s_*+q_*$, as required.
\end{proof}

%\section{Localization}\label{sec:local}

If $i:A\to A/I$ is a surjection of pointed monoids, there is an
exact functor $(A/I)\Sets\to A\Sets$. If $A$ is a pc monoid,
%so is $A/I$, and
we have a map $i_*: K'(A/I)\to K'(A)$.

\begin{thm}\label{A-A/s}
Let $A$ be a pc monoid. Suppose that $s\in A$ is such that the set $\{s^n\}_{n\geq 0}$
is a $2$--sided denominator set (that is, $s$-left fractions and $s$-right fractions coincide,
see \cite[II.A]{WK}).
%Then for each $s\in A$
%, and $s\in A$. %a nonzero-divisor.
%\edit{why a nzd?}
Then we have a fibration sequence
\[
K'(A/sA) \map{i_*} K'(A) \map{} K'(A[s^{-1}]).
\]
\end{thm}

\begin{proof}
Consider the category $\cT$ of (pc) $s$-torsion $A$--sets;
by Corollary \ref{cor-ACGW}, $\cT$
is an ACGW-category.
In fact, it is a sub ACGW-category of $A\Sets_\pc$.
Since every $s$-torsion $A$--set $X$ is finitely generated, $X$
has a finite filtration by the subsets $s^nX$. By Devissage,
$K'(A/sA) \simeq K(\cT)$.

We will apply the localization theorem \cite[8.5]{CZ} to
$\cT \subset A\Sets_\pc$.
The argument of \cite[8.6]{CZ} shows that
the category $A\Sets_\pc\backslash\cT$ is equivalent to $A[1/s]\Sets_\pc$,
and hence is a CGW category, so axiom (CGW) holds.
Axiom (W) holds: $\cT$ is ``m-negligible''
in $A\Sets_\pc$ because every pc $A$--set $N$ contains $s^n N$
(which is $s$-torsion-free for sufficiently large $n$), and if the kernel of $M\onto N$ is in $\cT$, then $s^n M\map{\sim} s^n N$ for $n$ large enough. Moreover, $\cT$ is
``e-well represented'' because, given pc $A$--sets $N_1$, $N_2$ and $V$, and $A$--set
maps
\[
	N_1\to N_1[s^{-1}] \map{\sim} V[s^{-1}]\; \text{and}\; N_2\to N_2[s^{-1}] \map{\sim} V[s^{-1}]
\]	
	 the pullback $N = N_1\times_{V[s^{-1}]} N_2$ is a pc $A$--set and $N[s^{-1}]\to V[s^{-1}]$ is an isomorphism.
%because every morphism in the category $A\Sets_\pc$
%can be factored into an admissible epimorphism followed by an admissible monomorphism,
%and the pullback of any morphism that
%becomes an isomorphism in $A[1/s]\Sets_\pc$ also becomes an isomorphism after localization.
%the quotient $N'$ of $N$ by its
%$s$-torsion $A$--subset has no $s$-torsion.
%is $s$-torsionfree.
%``e-negligible'' because the pullback of an $s$-torsion
%$A$--subset is an $s$-torsion $A$--set; see Definition 8.4 in \cite{CZ}.
 Finally, axiom (E) holds because every morphism $A\to B$ in
$A\Sets_\pc\backslash\cT$ is represented as
$A \onto A/(s\mathrm{-torsion}) \into s^nB$ for some $n$.
\end{proof}

The following corollary is an analogue of \cite[Ex.\,V.6.4]{WK}, \cite[Exercise in \S 6]{Q} and \cite{FarrellHsiang}.

\begin{cor}\label{local}
Let $A$ be a pc monoid and $\phi: A\to A$ an automorphism. Write $A\rtimes\N$ and $A\rtimes\Z$, respectively, for the corresponding semi-direct product monoids (in which $ta = \phi(a) t$). Write $i: A\to A\rtimes\N$ for the inclusion. Then there is a fibration sequence
\[ \xymatrix@C=3em{
	K'(A) \ar[r]^{i_*-i_*\phi_*~} &  K'(A\rtimes\N) \ar[r] & K'(A\rtimes \Z_+).
}
%K'(A) \map{i_* - i_* \phi_*} K'(A\rtimes\N) \smap{} K'(A\rtimes\Z).
\]
In particular (taking $\phi$ to be the identity), we get a fibration sequence
\[
K'(A) \map{0} K'(A\wedge\N) \smap{} K'(A\wedge\Z_+).
\]
Thus $K'(A\wedge\Z_+)\simeq K'(A\wedge\N)\vee \Omega^{-1}K'(A)$.
\end{cor}

\begin{proof}
Apply Theorem \ref{A-A/s} to $A\rtimes\N$ and $s = t$. The fact that
$K'(A) \to K'(A\rtimes\N)$ is $i_* - i_* \phi_*$ follows from additivity applied to
the characteristic exact sequence of an $A$--set:
$\phi_*(X)\wedge\N \overset{t}{\into} X\wedge\N \overset{i}{\onto} X$, where
%\edit{Check twist! Checked by Chuck}
the action of $A\rtimes\N$ on $\phi_*(X)\wedge\N$ is twisted by $\phi$.
\end{proof}

\section{Fundamental Theorem}\label{sec:FT}

We now apply the results in the previous section to prove
Theorems \ref{Fund.Thm} and \ref{thm:twisted}.

\begin{lem}\label{spectra}
Let $R$, $S$ and $T$ be spectra of finite type. Given an equivalence
$f:S\vee T\to R\vee S\vee T$, which is an equivalence on $S$, then
$f$ induces an equivalence $T\simeq T$ and $R\simeq0$.
\end{lem}

\begin{proof}
Localizing at $\Q$, a dimension count shows that the homotopy groups
of $R$ are finite abelian groups. Taking coefficients modulo a prime,
a cardinality argument shows that all homotopy groups of $R$ vanish.
Hence $R\simeq0$, and the result follows easily.
\end{proof}

\begin{lem}\label{split}
%Suppose $A$ is a pc monoid.
%$\phi:A\to A$ is an automorphism, and $A\rtimes \N$ and
%$A\rtimes \Z$ are the corresponding semidirect product monoids.
%\edit{works for the semidirect products. Does it though?}
The composition $K'(A) \map{f^*} K'(A\wedge\N) \map{j^*}K'(A\wedge\Z_+)$
splits, where $f^*$ and $j^*$ are induced from the base change
%\edit{is this lemma used anywhere?}
$X\mapsto X\wedge_A \N$ and localization functor $X\mapsto X[1/t]$.
\end{lem}

\begin{proof}
Since both $f^*$ and $j^*$ are exact, the composition is defined.
Consider the orbit functor $\gamma$ from $A\wedge\Z_+\Sets_\pc$ to $A\Sets_\pc$,
sending $X$ to $X/\sim$, where $x\sim y$ if $y=t^nx$ for some integer $n$.
It is easy to see that $\gamma$ is an exact functor, and that
$\gamma j^*f^*$ is the identity.
\end{proof}

\begin{thm}\label{K'N}
$K'(\N) \simeq \bbS$.
\end{thm}

\begin{proof}
Let $\cC$ denote the category of finite pc $\N$--sets;
it is an $ACGW$ category by Theorem \ref{ACGW}.
%The inclusion of $\cC$ into pc $\N$--sets is zero
%by Corollary \ref{local}.
Since each finite pc $\N$--set $X$ has a finite filtration by
$t^nX$, devissage \ref{thm:devissage} and Example \ref{BBQ} imply
that $K(\cC)\simeq K(\mathbf{Sets}_\fin)\simeq \bbS$.
%the additivity follows from  \cite[7.15]{CZ} because $\cC$ is
%\edit{needs ACGW}
%an $ACGW$ category.
Then the localization sequence \ref{local} %\cite[8.5]{CZ}
yields a fibration sequence
\[
K(\mathbf{Sets}_\fin) \map{0} K'(\N) \map{} K'(\Z_+).
\]
Hence $K'(\Z_+)\simeq K'(\N)\vee \Omega^{-1}K(\mathbf{Sets}_\fin)$.
By Example \ref{K(G)}, we see that
$K'(\Z_+)\simeq \bbS \vee \Omega^{-1}K(\mathbf{Sets}_\fin)$.
Since the inclusion (given by base extension) $\bbS \simeq K(\mathbf{Sets}_\fin)\to K'(\Z_+)$ factors through $K'(\N)$,
Lemma \ref{split} implies that there
is a spectrum $R$ such that $K'(\N)\simeq R\vee \bbS.$
By Lemma \ref{spectra},
$K'(\mathbf{Sets}_\fin)\simeq K'(\N)$.
\end{proof}

\begin{porism}\label{K'GN}
If $G$ is a group, then the proof of Theorem \ref{K'N} applies,
using Lemma \ref{lem:length}, to
show that $K'(G_+\wedge\N)\simeq K'(G_+)\simeq S^\infty(BG_+)$.
\end{porism}

\begin{thm}\label{GxZ}
If $A$ is a pc monoid of finite length with units $G$, then
\[
K'(A\wedge\Z_+)\cong K'(G\times\Z_+) \simeq S^\infty(BG_+)
\vee \Omega^{-1}S^\infty(BG_+).
  \]
In particular, if $G=\{1\}$ then $K'(A\wedge\Z_+)\cong K'(\Z_+)\simeq \bbS\vee \Omega^{-1}\bbS$.
\end{thm}

\begin{proof}
	$A\wedge\Z_+$ has finite length, and its units are $G\times\Z.$ By Lemma \ref{lem:length},
	we conclude that $K'(A\wedge\Z_+)\cong K'(G\times\Z_+).$
%Let $\cC$ denote the category of pc $A\wedge\Z_+$--sets
%on which $\m_A$ acts trivially; $\cC$ is equivalent to pc $(G\times\Z)_+$--sets.
%Every (pc) $A\wedge\Z_+$--set $X$ has  a finite filtration by $\m_A^iX$,
%and the quotients are in $\cC$. %(i.e., pc $\Z_+$--sets).
%Now apply devissage \ref{thm:devissage}; $K(\cC)\simeq K(A\wedge\Z_+\Sets_\pc))=K'(A\wedge\Z_+)$.
By Example \ref{K(G)}, this is $S^\infty B(G\times\Z)_+$, and
\[
 S^\infty B(G\times\Z)_+ \simeq
 S^\infty(BG_+)\vee \Omega^{-1} S^\infty(BG_+).\qedhere
\]
\end{proof}

%\begin{cor}
%If $G$ is a finitely generated abelian group
%\end{cor}

\begin{thm}\label{Fund.Thm}
If $A$ is a pc abelian monoid, then
\[K'(A)\simeq K'(A\wedge\N)\quad \textrm{and}\quad K'(A\wedge\Z_+)\simeq K'(A)\vee \Omega^{-1}K'(A).\]
\end{thm}

\begin{proof}
	It suffices to prove the first equivalence, as the
second equivalence follows from it, using Corollary \ref{local}.
If $A=G_+$ then $K'(G_+)\simeq K'(G_+\wedge\N)$ by Porism \ref{K'GN}.
%the proof of Theorem \ref{K'N} applies,
%replacing $Omega^\infty S^\infty$ with $Omega^\infty S^\infty(BG_+)$.
Inductively, suppose that the result holds for pc monoids with
$n$ generators over its group of units, and that $A$ is generated
by $s=s_0,...,s_n$, where $s$ is in $\m_A$.
Then the result holds for the monoids $A/sA$ and $A[1/s]$ (which are both pc by \cite[Prop. 9.1]{CHWW}),
so by naturality and Theorem \ref{A-A/s} (which applies as $A$ is abelian) we have a map of fibration sequences
whose outside maps are equivalences:
\begin{equation*}
\xymatrix@R=1.9em{
K'(A/sA) \ar[r]\ar[d]^{\simeq} & K'(A)\ar[r]\ar[d]& K'(A[1/s]) \ar[d]^{\simeq}\\
K'(A/sA\wedge\N) \ar[r] & K'(A\wedge\N)\ar[r]& K'(A[1/s]\wedge\N).
}\end{equation*}
By the 5--lemma, the middle map $K'(A)\to K'(A\wedge\N)$ is an equivalence.
%Now the second equivalence follows using Corollary \ref{local}.
\end{proof}

\begin{thm}\label{thm:twisted}
Let $G$ be a group,
$\phi: G\to G$ an automorphism, and $G_+\rtimes \N$ and $G_+\rtimes\Z$
the corresponding semidirect product monoids.
Then the base extension map $K'(G_+)\to K'(G_+\rtimes\N)$ is an equivalence
and we have a fibration sequence $K'(G_+) \map{1 - \phi_*} K'(G_+) \to K'(G_+\rtimes\Z).$ 	
\end{thm}

\begin{proof}
We consider the diagram
\begin{equation*}
\xymatrix@R=1.9em@C=3em{
	K'(G_+) \ar[r]^{1 - \phi_*}\ar[d]^{=} & K'(G_+)\ar[r]\ar[d]^{i_*}& K'(G_+\rtimes\Z) \ar[d]^{=}\\
	K'(G_+) \ar[r]^{i_* - i_* \phi_*} & K'(G_+\rtimes\N)\ar[r]& K'(G_+\rtimes\Z),
}
\end{equation*}
where $i_*$ is the base extension map along the inclusion $G_+\to G_+\rtimes\N$.
The bottom sequence is the localization sequence of Corollary \ref{local}, and the whole diagram commutes by inspection.
It therefore suffices to show that the top row is a cofibration sequence; it is canonically equivalent to the sequence of
suspension spectra of classifying spaces $S^\infty BG_+ \map{1 - \phi_*} S^\infty BG_+ \to S^\infty B(G\rtimes \Z)_+.$
This last sequence is a cofibration sequence by the mapping torus construction:
unstably, the homotopy coequalizer of $1$ and $\phi_*$ is computed by the
mapping torus $T = BG\times I/\{(x,0) \sim (\phi_*(x), 1)\}.$ This mapping torus
has $BG\times\R$ as a covering space with deck transformation group $\Z$.
It follows that $\pi_n(T) = 0$ for $n\geq 2$, and it is easy to verify
that $\pi_1(T) \cong G\rtimes \Z.$ Thus $T\simeq B(G\rtimes\Z)$, and stabilisation yields the desired
conclusion.
\end{proof}

\medskip

\section{Monoid schemes}\label{sec:MS}

The constructions of Section \ref{sec:K'} can be
generalized to define the $K'$-theory of noetherian partially cancellative monoid schemes.
For relevant definitions see \cite{CHWW}. In particular, given an abelian pc monoid $A$,
we have an affine pc monoid scheme $\mathrm{MSpec}(A)$; a general pc monoid scheme is a noetherian ringed space
locally isomorphic to such an affine pc monoid scheme. Recall as well from \cite[Defn.\,2.5]{CHWW}
that an equivariant closed subscheme
of a monoid scheme $X$ is one that is defined by a sheaf of ideals.

\begin{defn}
Let $X$ be a pc monoid scheme. A pc set on $X$ is a sheaf
that is locally of the form $\widetilde{F}$ for a finitely generated pc $A$--set $F$; compare \cite[3.3]{CLS}.
We write $X\Sets_\pc$ for the (quasi-exact) category of pc sets on $X$;
its admissible sequences are those that are locally
admissible in the sense of Example \ref{G-and-K'}.
We define $K'(X)$ to be the $K$-theory of $X\Sets_\pc.$
\end{defn}

\begin{subrem}
If $X = \mathrm{MSpec}(A)$, then there is a natural equivalence of quasi-exact categories $X\Sets_\pc \cong A\Sets_\pc$
(see \cite[Cor.\,3.13]{CLS}), and therefore $K'(X)\simeq K'(A). $
\end{subrem}

\begin{subrem}
The $K'$-theory of pc monoid schemes is contravariantly functorial for flat morphisms,
and covariantly functorial for equivariant closed immersions. This is immediate from Remark \ref{K'-functor}.
\end{subrem}

As pushouts and pullbacks of sheaves, and admissible sequences,
are all detected locally,
the following result is an immediate consequence of Theorem \ref{ACGW}.

\begin{thm}\label{ACGW-schemes}
Let $X$ be a pc monoid scheme. Then any quasi-exact subcategory of $X\Sets_\pc$ that is closed under pushouts
and pullbacks is an ACGW-category.	
\end{thm}

Here is the analogue of Quillen's localisation theorem \cite[7.3.2]{Q}.

\begin{thm}\label{thm:X-U}
Let $X$ be a pc monoid scheme and $Z\map{i} X$ an equivariant closed subscheme
with open complement $U\map{j} X$. Then there is a fibration sequence of spectra
\[
K'(Z)\map{i_*} K'(X) \map{j^*} K'(U).
\]
\end{thm}

\begin{proof}
Let $J$ be the ideal sheaf defining $Z$. The proof of Theorem \ref{A-A/s}
applies {\it mutatis mutandis}, with $\cT$ the subcategory of $X\Sets_\pc$
consisting of those sheaves supported on $Z$, or equivalently,
the subcategory of $J$-torsion sheaves. By devissage, $K(\cT)$
is equivalent to $K'(Z)$. The localized (bi-)category $X\Sets_\pc\backslash\cT$
is equivalent to $U\Sets_\pc$; in particular, axiom (CGW) of \cite[Theorem 8.5]{CZ} holds.
The remaining axioms are checked as in the proof of Theorem \ref{A-A/s}.
\end{proof}	

\begin{subrem}\label{MV}
The localisation sequence of the theorem is natural for flat pullbacks $X'\to X$.
In particular, if $V\subset X$ is an open subscheme containing $Z$, then
\begin{equation*}
\xymatrix@R=1.9em@C=3em{
	K'(Z) \ar[r]\ar[d]^= & K'(X)\ar[r]\ar[d]& K'(U) \ar[d]\\
	K'(Z) \ar[r] & K'(V)\ar[r]& K'(U\cup V).
}
\end{equation*}
is a map of fibration sequences of spectra.
\end{subrem}

The global version of the Fundamental Theorem is now easily derived.
Recall that $\mathbb{A}^1$ denotes $\mathrm{MSpec}(\N)$.

\begin{thm}\label{Fund.Thm:X}
Let $X$ be a pc monoid scheme. Then the pullback map $K'(X)\to K'(X\times\mathbb{A}^1)$ is an equivalence,
and
\[K'(X\wedge\Z_+)\simeq K'(X)\vee\Omega^{-1}K'(X).\]
\end{thm}

\begin{proof}
Note that $\mathrm{MSpec}(A)\times\mathbb{A}^1 = \mathrm{MSpec}(A\wedge\N)$.
By Theorem \ref{Fund.Thm} the assertion holds for affine pc monoid schemes,
and in particular for pc monoid schemes of dimension $0$.
The general case is now proved using induction on the dimension of $X$ using Theorem \ref{thm:X-U}.
\end{proof}

\begin{cor}\label{homotopy}
Let $\xi:E\to X$ be a vector bundle over a pc monoid scheme. Then $\xi^*: K'(X)\to K'(E)$ is an equivalence.
\end{cor}

\begin{proof}
For a trivial vector bundle, this follows from the theorem by induction on the rank.
Now the general case follows from localization \ref{thm:X-U},
applied to an open cover of $X$ trivializing $E$.
\end{proof}	

\begin{thm}
Let $X$ be a pc monoid scheme. Then there is a natural equivalence $K'(X\times\mathbb{P}^1)\simeq K'(X)\vee K'(X)$.
\end{thm}

\begin{proof}
We have an open covering $X\times\mathbb{P}^1 = U_1\cup U_2$,
with each $U_i \cong X\times\mathbb{A}^1$ and $U_1\cap U_2\cong X\wedge\Z_+$.
The complement of each $U_i$ is isomorphic to $X$.
Applying Remark \ref{MV} and the Fundamental Theorem \ref{Fund.Thm:X},
we obtain a homotopy cocartesian square of spectra
\begin{equation*}
\xymatrix{ %}@R=1.9em@C=3em{
	 K'(X\times\mathbb{P}^1)\ar[r]\ar[d]& K'(X) \ar[d]\\
	 K'(X)\ar[r]& K'(X)\vee \Omega^{-1}K'(X).
}
\end{equation*}
The right vertical and bottom map are isomorphisms onto the first summand; the assertion follows.
\end{proof}	

\begin{thm}
	Let $X$ be a pc monoid scheme, and $n$ a non-negative integer. Then there is a natural equivalence
	\[K'(X\times\mathbb{P}^n)\simeq K'(X)^{\vee n+1}.\]
	In particular, $K'(\mathbb{P}^n) \simeq \bbS \vee \bbS \vee \cdots \vee \bbS$  ($n+1$ copies of $\bbS$).
\end{thm}

\begin{proof}
We proceed by induction on $n$, the case $n = 0$ being tautological.
Suppose $n > 0$ and we have proved the assertion for $n-1$. Applying Theorem \ref{thm:X-U}
to the equivariant closed immersion $X\times\mathbb{P}^{n-1}\to X\times\mathbb{P}^n$ with open
complement $U \cong X\times\mathbb{A}^n$ we obtain a fibration sequence of spectra
\[
K'(X\times\mathbb{P}^{n-1}) \to K'(X\times\mathbb{P}^{n}) \to K'(X\times\mathbb{A}^{n})\simeq K'(X).
\]
Because $K'(X)\to K'(X\times\mathbb P^n) \to K'(X\times\mathbb A^n)$ is an equivalence  by Corollary \ref{homotopy},
this sequence is split; the assertion for $n$ follows.
\end{proof}

\begin{subrem}
In particular, $K'_0(\mathbb{P}^1)\cong \Z\oplus\Z\cong \Z\times\mathrm{Pic}(\mathbb{P}^1)$,
and the isomorphism is given by the rank and determinant. See \cite{FW} for the computation of
the Picard group of monoid schemes.
\end{subrem}

%One new ingredient here is
%the Picard group, as investigated in \cite{FW}. See \cite[Section 4]{Sz} for
%description of sheaves on projective line.

\medskip
\subsection*{Acknowledgements}
The authors would like to thank Inna Zakharevich for her help in
understanding CGW categories in \cite{CZ}.

\medskip\goodbreak

\end{document}